\documentclass[leqno,12pt]{amsart} 

\usepackage{amssymb,epsfig,verbatim} 

\theoremstyle{plain}
\newtheorem{theorem}{Theorem}[section]
\newtheorem{proposition}[theorem]{Proposition}
\newtheorem{corollary}[theorem]{Corollary}
\newtheorem{lemma}[theorem]{Lemma}

\theoremstyle{remark}
\newtheorem{definition}[theorem]{Definition}

\newtheorem{remark}[theorem]{Remark}

\newcommand{\s}{{\mathbb S}}
\renewcommand{\r}{{\mathbb R}}
\newcommand{\h}{{\mathbb H}}

\renewcommand{\c}{{\mathbb C}}

\DeclareMathOperator{\trace}{trace}
\DeclareMathOperator{\grad}{grad}
\DeclareMathOperator{\ricci}{Ricci}

\DeclareMathOperator{\cst}{constant}
\renewcommand{\sl}{{\widetilde{\mathrm{SL}_{2}(\r)}}}
\newcommand{\sol}{{\mathrm{Sol}}}

\begin{document}
\title{Biminimal Immersions}

\author{E. Loubeau}
\address{D{\'e}partement de Math{\'e}matiques \\
Laboratoire CNRS UMR 6205 \\
Universit{\'e} de Bretagne Occidentale \\
6, avenue Victor Le Gorgeu \\
CS 93837, 29238 Brest Cedex 3\\
France}
\email{loubeau@univ-brest.fr}

\author{S. Montaldo}
\address{Dipartimento di Matematica e Informatica \\
Via Ospedale 72 \\
09124 Cagliari \\
Italy}
\email{montaldo@unica.it}

\dedicatory{Dedicated to Professor Renzo Caddeo on his 60th
birthday}
\thanks{Work partially supported by GNSAGA(Italy) and University of Brest.}

\date{\today}

\maketitle
\begin{abstract}
We study biminimal immersions, that is immersions which are critical points of the bienergy
for normal variations with fixed energy.  We give a geometrical description of the
Euler-Lagrange equation associated to biminimal immersions for: i) biminimal curves in a Riemannian
manifold, with particular care to the case of curves in a space form
ii) isometric immersions of
codimension one in a Riemannian manifold, in particular  for surfaces of a three-dimensional manifold. We describe two methods to construct families of biminimal surfaces using both Riemannian and horizontally homothetic submersions.
\end{abstract}
\section{Introduction}

Many stimulating problems in mathematics owe their existence to variational
formulations of physical phenomena.
In differential geometry, harmonic maps, candidate minimisers of the Dirichlet energy, can
be described as constraining a rubber sheet to fit on a marble manifold in a
position of elastica equilibrium, i.e. without tension~\cite{E-L1}.
However, when this scheme falls through, and it can, as corroborated by the case of the
two-torus and the two-sphere~\cite{E-W}, a best map will minimise this failure, measured by the total tension, called bienergy.
In the more geometrically meaningful context of immersions, the fact that the 
tension field is normal to the image submanifold, suggests that the most effective
deformations must be sought in the normal direction.
\newline Then, two approaches to this optimization are available: The first one (the free state) consists in finding (normal) extrema of the bienergy, with complete disregard for the behaviour of the energy. In the second, in order to avoid paying too great a price for a smaller tension, a constancy condition on the energy level is imposed.
\newline In more intuitive terms, and even though we never consider the associated flows, these points of view correspond, at least in the more favorable situations, to reducing the overall tension of a surface, with or without controlling the mean curvature.
\newline However different, both outlooks are unified into a single mathematical description, amounting to a Lagrange multiplier interpretation.

This considerations lead to the following definitions.

\begin{definition}
A map $\phi : (M,g) \to (N,h)$ between Riemannian 
manifolds is called {\em biharmonic} if it is 
a critical point, for all variations,
of the bienergy functional: $$
E_{2}(\phi) = \frac{1}{2} \int_{M} |
\tau(\phi)|^2 \, v_{g} ,$$ where
$\tau(\phi) = \trace{\nabla d\phi}$ is
the tension field, vanishing for critical points of the Dirichlet energy (i.e. harmonic maps):
$$E(\phi) = \frac{1}{2} \int_{M} |d\phi|^2 \, v_{g} .$$ 
In case of non-compact domain, these two definitions should be understood as for all compact subsets. 
\end{definition}

The Euler-Lagrange operator attached to biharmonicity, called
the {\em bitension field} and  computed by Jiang in~\cite{Jiang}, is:
$$ 
\tau_{2}(\phi) = -\left(\Delta^{\phi}
\tau(\phi) - \trace R^{N}(d\phi
,\tau(\phi) ) d\phi \right),
$$ 
and vanishes if and only if the map $\phi$ is biharmonic.

We are now ready to define the main object of this paper.

\begin{definition}
An immersion $\phi : (M^{m},g) \to (N^{n},h)$ ($m \leq n$) between 
Riemannian manifolds, 
or its image, is called {\em
biminimal} if it is a critical point of
the bienergy functional $E_2$ for
variations normal to the image
$\phi(M)\subset N$, with fixed energy. Equivalently, there exists a
constant $\lambda\in\r$ such that $\phi$ is a critical point of
the $\lambda$-bienergy
$$
E_{2,\lambda}(\phi) = E_{2}(\phi) + \lambda E(\phi)
$$
for any  smooth variation of the map $\phi_{t} : ]-\epsilon ,
+\epsilon[ \times M \to N$, $\phi_{0}=\phi$, such that $V=\frac{d
\phi_{t}}{dt}|_{t=0}$ is normal to $\phi(M)$. 
\end{definition}

\begin{remark}
The functional $E_{2,\lambda}$ has been on the mathematical scene since the early seventies 
(see~\cite{Eliasson}), where its critical points, for all possible variations, are studied. In particular, it is shown to satisfy Condition~(C) when the domain has dimension two or three and the target is non-positively curved, ensuring the existence of minimisers in each homotopy class. However, L.~Lemaire, in~\cite{Lemaire}, constructs counter-examples when no condition is imposed on the curvature.
\end{remark}

Using the Euler-Lagrange equations for harmonic and biharmonic maps, 
we see that an immersion is biminimal if 
$$
[\tau_{2,\lambda}]^{\perp}=[\tau_{2}]^{\perp}- \lambda [\tau]^{\perp}=0,
$$
for some value of $\lambda \in \r$, where $[,]^{\perp}$ denotes the normal component of $[,]$.

We call an immersion {\it free} biminimal if it is
biminimal for $\lambda=0$.

In the instance of an isometric  immersion $\phi : M \to N$,
the biminimal condition is 
\begin{equation}\label{eq:bi-is}
[\Delta^{\phi}
{\mathbf H} - \trace R^{N}(d\phi
,{\mathbf H} ) d\phi]^{\perp}+\lambda {\mathbf H} = 0.
\end{equation}

\noindent Note that this variational principle
is close to the Willmore problem, the disparity being that we do not
vary through isometric immersions.

While it is obvious that biharmonic immersions 
are biminimal, we will see in the following sections that the two notions are well distinct. For example, we construct families of biminimal surfaces
in three-dimensional space forms of non-positive constant sectional curvature where biharmonic surfaces do not exist \cite{CBYIS,CMO2}.
In the same vein,  we construct families of biminimal
surfaces in almost all three-dimensional geometries of Thurston.

\noindent{\it Notation}.
We shall place ourselves in the $C^{\infty}$ category, 
i.e. manifolds, metrics, 
connections, maps will be  assumed to be smooth.
By $(M^m,g)$ we shall mean a connected manifold,
 of dimension $m$, without boundary, endowed with a 
 Riemannian metric $g$. 
We shall denote by $\nabla$ the Levi-Civita connection 
on $(M,g)$. For vector fields $X,Y,Z$ on $M$ we define the 
Riemann curvature operator by $R(X,Y)Z=\nabla_{[X,Y]}Z-[\nabla_{X},\nabla_{Y}]Z$.
For the Laplacian we shall use $\Delta(f)={\rm div}\grad f$ for
functions $f\in C^{\infty}(M)$ and
$\Delta^{\phi}W=-\trace\left(\nabla^{\phi}\right)^2\,W$ for
 sections along a map
$\phi:M\to N$.

\section{Biminimal curves}

Our quest for examples of  biminimal 
immersions starts with curves.

Let $\gamma : I\subset {\r}\to (M^{m},g)$ be a curve parametrised
by arc-length in a Riemannian manifold $(M^{m},g)$,
that is $\gamma$ is an isometric immersion. Before computing the
bitension field of  $\gamma$, we recall the definition of Frenet
frames.

\begin{definition}[See, for example, \cite{Laugwitz}] \label{def2.1}
The Frenet frame $\{B_{i}\}_{i=1,\dots,m}$ associated to a curve 
$\gamma : I\subset {\r}\to (M^{m},g)$ is the orthonormalisation of
the $(m+1)$-uple  $\{ \nabla_{\frac{\partial}{\partial
t}}^{(k)} d\gamma
(\frac{\partial}{\partial t})
\}_{k=0,\dots,m}$, described by:
\begin{align*} B_{1}&=d\gamma
(\frac{\partial}{\partial t}) ,  \\
\nabla_{\frac{\partial}{\partial
t}}^{\gamma} B_{1} &= k_{1} B_{2} , \\
\nabla_{\frac{\partial}{\partial
t}}^{\gamma} B_{i} &= - k_{i-1} B_{i-1}
+ k_{i}B_{i+1} , \quad \forall i =
2,\dots,m-1 , \\
\nabla_{\frac{\partial}{\partial
t}}^{\gamma} B_{m} &= - k_{m-1} B_{m-1}
, \end{align*} 
where the functions $\{k_{1}=k>0,k_{2}=\tau,k_{3},\ldots,k_{m-1}\}$
 are called the curvatures of $\gamma$. Note that
$B_{1}=T=\gamma'$ is the unit tangent vector field to the curve.
\newline In
the instance of a curve $\gamma$ on a surface ($m=2$), the Frenet
frame reduces to the couple
$\{T,N\}$, $T$ being
the unit tangent vector field along $\gamma$ and $N$ a normal vector
field along $\gamma$ such that $\{T,N\}$ is a positive basis,
 while $k_{1}=k$ is  the signed curvature of
$\gamma$.
 \end{definition}

Biminimal curves in a Riemannian manifold are characterised by:
 
\begin{proposition}\label{prop1}
Let $\gamma : I \subset {\r} \to (M^{m},g)$ ($m\geq 2$) 
be an isometric
curve from an open interval of $\r$ into a Riemannian 
manifold $(M,g)$. 
Then $\gamma$ is biminimal if and only
if there exists a real number $\lambda$ such that: 
\begin{equation}\label{sys1}
\left\lbrace
\begin{array}{l}
k_{1}'' - k_{1}^3 - k_{1}k_{2}^2 + k_{1}
g(R(B_{1},B_{2})B_{1},B_{2})-\lambda k_{1} = 0 , 
\\ (k_{1}^{2}k_{2})' + k_{1}^2
g(R(B_{1},B_{2})B_{1},B_{3}) = 0 , \\ 
k_{1}k_{3} + k_{1}
g(R(B_{1},B_{2})B_{1},B_{4}) = 0 , \\
k_{1} g(R(B_{1},B_{2})B_{1},B_{j}) = 0,\quad j=5,\ldots,m,
\end{array}
\right.
\end{equation}
where $R$ is the curvature tensor of $(M,g)$ and $\{B_{i}\}_{i=1,\dots,m}$ 
the Frenet frame of $\gamma$.
\end{proposition}

\begin{proof}
With respect to its Frenet frame, the tension field of $\gamma$ is:
\begin{align*}
\tau(\gamma) &= \trace{\nabla d\gamma} = \nabla_{\frac{\partial}{\partial t}}^{\gamma} (d\gamma
(\frac{\partial}{\partial t})) - d\gamma (\nabla_{\frac{\partial}{\partial t}}\frac{\partial}{\partial t}) = \nabla_{\frac{\partial}{\partial t}}^{\gamma}
B_{1} = k_{1} B_{2} 
\end{align*}
and its bitension field:
\begin{align*}
-\tau_{2}(\gamma) =& -\nabla_{\frac{\partial}{\partial
t}}^{\gamma}\nabla_{\frac{\partial}{\partial t}}^{\gamma}
(\tau(\gamma)) + \nabla_{\nabla_{\frac{\partial}{\partial
t}}\frac{\partial}{\partial t}}^{\gamma} (\tau(\gamma)) - 
R(d\gamma (\frac{\partial}{\partial t}),\tau(\gamma) )d\gamma
(\frac{\partial}{\partial t})\\
=& -\nabla_{\frac{\partial}{\partial
t}}^{\gamma}\nabla_{\frac{\partial}{\partial t}}^{\gamma} (k_{1}
B_{2}) -R(B_{1},k_{1} B_{2})B_{1} \\ =&
-\nabla_{\frac{\partial}{\partial t}}^{\gamma} (k_{1}' B_{2} -
k_{1}^{2}B_{1} + k_{1}k_{2} B_{3}) - k_{1} R(B_{1},B_{2})B_{1} \\ =&
-(k_{1}'' - k_{1}^3 - k_{1}k_{2}^{2}) B_{2} + 3 k_{1}k_{1}' B_{1} -
(k_{1}'k_{2} +  (k_{1}k_{2})')B_{3} \\
& - k_{1}k_{3} B_{4} - k_{1}
R(B_{1},B_{2})B_{1}  .  
\end{align*} 
The vanishing of the normal
components  yields System~\eqref{sys1}. 
\end{proof}

\begin{remark}
Asking for a free biharmonic curve $\gamma$ to be biharmonic requires the supplementary condition
$[\tau_{2}(\gamma)]^{B_{1}} =0$, equivalent to $k_{1}k_{1}'=0$, that
is, either $k_{1}$ is constant or  $\gamma$ is a geodesic
($k_{1}=0$).  
\end{remark}

If the target manifold is a surface or a three-dimensional
Riemannian manifold with constant sectional curvature, Equations
\eqref{sys1} are more manageable as shown in the following:

\begin{corollary}\label{cor2.1}\quad\\
\begin{itemize}
\item[i)] An isometric curve $\gamma$ on a surface of Gaussian
curvature $G$  is biminimal if and only if its signed curvature $k$
satisfies the ordinary differential equation:
\begin{equation}\label{eq1} k'' - k^3 + kG -\lambda k= 0,
\end{equation}
for some $\lambda\in\r$.
\item[ii)]
An isometric curve $\gamma$ on a Riemannian three-manifold
 of constant sectional curvature $c$ 
is biminimal if and only if its curvature $k$ and torsion $\tau$
fulfill the system: 
\begin{equation}\label{eq2}
\left\{
\begin{array}{l}
k'' - k^3 -k\tau^2 + kc -\lambda k= 0 \\
k^2 \tau= \cst\,,
\end{array}
\right.
\end{equation}
for some $\lambda\in\r$.
\end{itemize}
\end{corollary}

\begin{proof} For i).
The two-dimensional Frenet frame of $\gamma$ consists 
only of $T$ and $N$, 
and the curve is biminimal, with respect to $\lambda$, if and only if :
$$ k'' - k^3 + k g(R(T,N)T,N) - \lambda k= 0 ,$$
but since $g(R(T,N)T,N) = G$, we obtain~\eqref{eq1}.

\noindent
In dimension three, the Frenet frame of $\gamma$ is
$\{T,N=B_{2},B=B_{3}\}$, and the conditions of 
Proposition~\ref{prop1}
become: \begin{equation*}
\left\lbrace
\begin{array}{l}
k'' - k^3 - k\tau^{2} + k g(R(T,N)T,N) - \lambda k= 0 , \\
(k^{2}\tau)' + k^2 g(R(T,N)T,B) = 0 .
\end{array}
\right.
\end{equation*}
The constant sectional curvature of the target means that
$g(R(T,N)T,N) = c$ and $g(R(T,N)T,B) =0$.
\end{proof}

>From Corollary~\ref{cor2.1}, if $\gamma$ is an isometric curve in a
Riemannian manifold $M^n$ of constant sectional curvature $c$ and
dimension $2$ or $3$, then the curvature of $\gamma$ (the signed
curvature when $n=2$), satisfies the equation:
\begin{equation}\label{eq3} 
k'' - k^3 -\frac{\alpha^2}{k^3}+ k\beta = 0,
\end{equation}
where $\alpha=k^2\tau$ and $\beta=c-\lambda$.

\begin{figure}[htb]
\epsfxsize=1.7in
\centerline{
\leavevmode
\epsffile{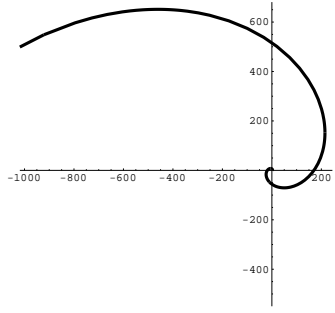}}
\caption{\label{figure1} Free biminimal curve in $\r^2$ of logarithmic
type.} \end{figure}

\noindent Multiplying Equation~\eqref{eq3} by $2k'$ and integrating, we
obtain: 
$$ 
(k')^2 - k^{4}/2 + \frac{\alpha^2}{k^2}+\beta k^2 = A ,\quad \; A\in\r,
$$
and setting $u = k^2$ yields:
$$ 
(u')^2 - 2 u^3 +4 \alpha^2+ 4  \beta u^2 = 4 Au.
$$
Since this equation is of the form $(u')^2= P(u)$, $P$ being a
polynomial of degree three, it can be solved by
standard techniques in terms of elliptic functions (see, for example
\cite{Davis}). In a forthcoming paper we shall give an accurate
description of the solutions of Equation~\eqref{eq3}. Here we just
point out that if  $M$ is the flat $\r^2$, then Equation~\eqref{eq3}, for
free biminimal curves,
reduces to  $$
k'' - k^3 = 0,
$$
a solution of which can be expressed in terms of elementary
functions, that is $k(s)=\sqrt{2}/s$, where $s$ is the arc-length.
Now using the standard formula to integrate a curve of known signed
curvature, we find that, up to isometries of $\r^2$, this free biminimal
curve is given by 
$$
\gamma(s)=s/3\left(\cos(\sqrt{2}\log s )+\sqrt{2}\sin(\sqrt{2}\log s ),
-\sqrt{2} \cos(\sqrt{2}\log s)+\sin(\sqrt{2}\log s)\right).
$$
This is the standard parametrisation by arc-length of the logarithmic
spiral plotted in Figure~\ref{figure1}.

\subsection{Biminimal curves via conformal changes of the metric}

On a Riemannian manifold $(M,g)$, any representative of the conformal class $[g]$  can be expressed as $\overline{g} = e^{2f} g$, $f\in
C^{\infty}(M)$, and the Levi-Civita connections are related, for
$X,Y \in C(TM)$, by (cf. \cite{Besse}): 
\begin{equation}\label{eq-conf-change}
\overline{\nabla}_{X}Y = \nabla_{X}Y + X(f)Y +
Y(f)X - g(X,Y)\grad{f} .
\end{equation}
 Observe that
a geodesic $\gamma$ on $(M,g)$ will not remain geodesic after a
conformal change of metric, unless the conformal factor is
constant since: $$ \overline{\nabla}_{\dot{\gamma}}\dot{\gamma} =
\nabla_{\dot{\gamma}}\dot{\gamma} + 2 \dot{\gamma}(f)\dot{\gamma} -
|\dot{\gamma}|^2 \grad{f}.$$

The following theorem gives a tool to construct free biminimal
curves.

\begin{theorem}\label{thm1}
Let $(M^{m},g)$ be a Riemannian manifold. Fix a point $p\in M^{m}$ 
and choose a 
function $f$ depending only on the geodesic distance from $p$. Then any
geodesic on $(M,g)$ going through $p$ will be a free biminimal curve on
$(M^{m},\overline{g}= e^{2f}g)$. \end{theorem}

\begin{proof}
Let $\gamma$ be a geodesic of $(M^{m},g)$ and let
$\{B_{i}\}_{i=1,\dots,m}$ be the associated Frenet frame  (cf. Definition~\ref{def2.1}). Since the function $f$
depends only on the geodesic distance from $p$, i.e. the
$B_{1}$-direction, $B_{i}f=0, \, \forall i=2,\dots,m$. Since
$\gamma$ is a geodesic on $(M,g)$, then:
$$
\nabla_{\frac{\partial}{\partial t}}^{\gamma} B_{1} = 0 ,
$$ and the
tension field of $\gamma$ with respect to the metric $\overline{g}=
e^{2f}g$ is: 
\begin{align*} \overline{\tau}(\gamma) = 
\overline{\nabla}_{\frac{\partial}{\partial t}}^{\gamma}
d\gamma(\frac{\partial}{\partial t}) =
\overline{\nabla}_{\frac{\partial}{\partial t}}^{\gamma} B_{1} =
\nabla_{\frac{\partial}{\partial t}}^{\gamma} B_{1} + 2 B_{1}(f)B_{1}
-\grad{f} , 
\end{align*} 
where in this last equality we have used \eqref{eq-conf-change}. Besides: 
$$
\grad{f} =
B_{1}(f)B_{1} + \sum_{i=2}^{m} B_{i}(f) B_{i} = B_{1}(f)B_{1},
$$  
thus
$\overline{\tau}(\gamma) = B_{1}(f)B_{1}$. 
\newline Still with respect
to $\overline{g}$, the bitension field of $\gamma$ is: 
\begin{align*}
\overline{\tau}_{2}(\gamma) &= - \Delta^{\gamma} \overline{\tau}(\gamma)
+ \trace{\overline{R}(d\gamma, \overline{\tau}(\gamma))d\gamma}
\notag\\ &= \overline{\nabla}_{\frac{\partial}{\partial
t}}^{\gamma}\overline{\nabla}_{\frac{\partial}{\partial t}}^{\gamma}
(B_{1}(f)B_{1})  -
\overline{\nabla}_{\nabla_{\frac{\partial}{\partial
t}}\frac{\partial}{\partial t}}^{\gamma} (B_{1}(f)B_{1}) +
\overline{R}(B_{1},B_{1}(f)B_{1})B_{1} \notag\\ &=
\overline{\nabla}_{\frac{\partial}{\partial
t}}^{\gamma}\overline{\nabla}_{\frac{\partial}{\partial t}}^{\gamma}
(B_{1}(f)B_{1}) \notag\\ &=
\overline{\nabla}_{\frac{\partial}{\partial t}}^{\gamma}
(B_{1}B_{1}(f)B_{1} + B_{1}(f)^{2} B_{1}) \notag\\ &=
(B_{1}B_{1}B_{1}(f))B_{1} + B_{1}B_{1}(f)B_{1}(f)B_{1} +
2B_{1}B_{1}(f)B_{1}(f)B_{1} + (B_{1}(f))^3 B_{1} \notag\\ &= [
B_{1}B_{1}B_{1}(f) + 3 B_{1}B_{1}(f)B_{1}(f) + (B_{1}(f))^3] B_{1}. 
\end{align*} 
So $\overline{\tau}^{2}(\gamma)$ has no
normal component and $\gamma$ is free biminimal on  $(M^{m},\overline{g}=
e^{2f}g)$. 
\end{proof}


\begin{corollary}
Let $r$ be the geodesic distance from a point $p\in (M,g)$, and
$f(r)=\ln(ar^2+br+c),\;a,b,c\in\r$. Then a geodesic on $(M,g)$
through $p$ becomes a biharmonic  map on $(M,\overline{g}=e^{2f}g)$.  
\end{corollary}

\begin{proof}
>From the proof of Theorem~\ref{thm1}, a geodesic on $(M,g)$
through $p$ becomes a biharmonic  map on $(M,\overline{g})$ if $f$ is
a solution of the ordinary differential equation:
$$
 f'''(r)+3f''(r)f'(r)+f'(r)^3=0.
$$
To solve this equation, put $y=f'$ to obtain
\begin{equation} \label{eq-conf-bihar2}
y''+3y'y+y^3=0.
\end{equation}
Then, using the transformation $y=x'/x$,
Equation~\eqref{eq-conf-bihar2}  becomes $x'''/x=0$ which has the
solution $x(r)=\bar{a}r^2+\bar{b}r+\bar{c},\;\bar{a},\bar{b},\bar{c}\in\r$. Finally, from $f(r)=\ln(d\;x(r)),\;d\in\r$, 
we find the desired $f$.
\end{proof}

\begin{figure}[htb]
\epsfxsize=2.6in
\centerline{
\leavevmode
\epsffile{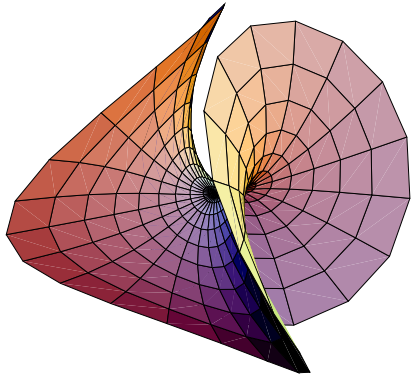}}
\caption{\label{enneper}The radial curves from the origin 
of the Enneper surface are biharmonic.} \end{figure}

\noindent As an example, one
can take $(M,g)=({\r}^{2},g=dx^2+dy^2)$, and $f(r)= \ln{(r^{2}
+1)}$,  where $r=\sqrt{x^2+y^2}$ is the distance from the origin.
 Thus any straight line on the flat
$\r^2$ turns into a biharmonic curve on $(\r^2,\bar{g}=(r^{2}
+1)^2(dx^2+dy^2))$ which is the metric, in local isothermal
coordinates, of the Enneper minimal surface. Figure~\ref{enneper}
is a plot of the Enneper surface in polar coordinates, so that radial curves on the picture are biharmonic.  

\section{Codimension-one biminimal submanifolds}

Let $\phi : M^n \to N^{n+1}$ be an isometric immersion 
 of codimension-one.
We denote by $B$ the second fundamental form of $\phi$,
by ${\mathbf N}$ a unit normal vector field to $\phi(M)\subset N$ and
by  ${\mathbf H}=H {\mathbf N}$ the
mean curvature vector field of $\phi$ ($H$ the mean curvature
function).

\noindent Then we have
\begin{proposition}\label{pro-B-ricci}
Let $\phi : M^n \to N^{n+1}$ be an isometric immersion 
of codimension one and ${\mathbf H} = H{\mathbf N}$ its mean
curvature vector. Then $\phi$ is
biminimal if and only if: 
\begin{equation}\label{eq5} 
\Delta H =
(|B|^2 - \ricci({\mathbf N})+\lambda) H. 
\end{equation}
for some value of $\lambda$ in $\r$.
\end{proposition}

\begin{proof}
In a local orthonormal
frame $\{e_{i}\}_{i=1,\dots,n}$ on $M$, the tension field of $\phi$ is
$\tau(\phi)= n H {\mathbf N}$ and its bitension field is:
\begin{align*} 
-\tau_{2}(\phi) =& -\sum_{i=1}^{n}\left[
\nabla^{\phi}_{e_{i}} \nabla^{\phi}_{e_{i}} (nH{\mathbf N}) +
\nabla^{\phi}_{\nabla_{e_{i}}e_{i}} (nH{\mathbf N})  -
  R^{N^{n+1}} (d\phi(e_{i}), nH{\mathbf N})d\phi(e_{i})\right]\\ 
=& n \sum_{i=1}^{n} \left[ -\nabla^{\phi}_{e_{i}}\left( e_{i}(H)
{\mathbf N} + H \nabla^{\phi}_{e_{i}}{\mathbf N} \right) 
+(\nabla_{e_{i}}e_{i})(H) {\mathbf N} +
H\nabla^{\phi}_{\nabla_{e_{i}}e_{i}} {\mathbf N}  - \right. \\
& \left.- H R^{N^{n+1}}(d\phi(e_{i}), {\mathbf N})d\phi(e_{i})
\right]\\  
=& n \sum_{i=1}^{n} \left[ -e_{i}e_{i}(H) {\mathbf N} - 2 e_{i}(H)
\nabla^{\phi}_{e_{i}}{\mathbf N} - H \nabla^{\phi}_{e_{i}}
\nabla^{\phi}_{e_{i}} {\mathbf N} +(\nabla_{e_{i}}e_{i})(H) {\mathbf
N} +\right.\\
& \left.+ H\nabla^{\phi}_{\nabla_{e_{i}}e_{i}} {\mathbf N}\right]
  - nH \sum_{i=1}^{n} R^{N^{n+1}} (d\phi(e_{i}), {\mathbf
N})d\phi(e_{i}) \\ 
=& n (\Delta H) {\mathbf N} - 2n \sum_{i=1}^{n}
e_{i}(H) \nabla^{\phi}_{e_{i}}{\mathbf N} + n H \Delta^{\phi}
{\mathbf N} - 
 nH\sum_{i=1}^{n} R^{N^{n+1}} (d\phi(e_{i}), {\mathbf
N})d\phi(e_{i}). 
\end{align*}  
But  
\begin{itemize} 
\item[i)] $\langle
\nabla^{\phi}_{e_{i}}{\mathbf N}, {\mathbf N}\rangle = \frac{1}{2}
e_{i}\langle{\mathbf N},{\mathbf N}\rangle =0$; 
\item[ii)] $\langle\sum_{i=1}^{n}
R^{N^{n+1}} (d\phi(e_{i}), {\mathbf N})d\phi(e_{i}), {\mathbf N}\rangle=
\ricci({\mathbf N})$.
\end{itemize}
For $\langle \Delta^{\phi}{\mathbf N},{\mathbf N}\rangle$, first we have:
$$\langle \Delta^{\phi}{\mathbf N},
{\mathbf N}\rangle = \sum_{i=1}^{n} \langle
-\nabla^{\phi}_{e_{i}}\nabla^{\phi}_{e_{i}}{\mathbf
N}+\nabla^{\phi}_{\nabla_{e_{i}}e_{i}}{\mathbf N},{\mathbf N}
\rangle=\sum_{i=1}^{n} \langle \nabla^{\phi}_{e_{i}}{\mathbf N},
\nabla^{\phi}_{e_{i}}{\mathbf N} \rangle.
$$ 
Then, if $B$ is the second
fundamental form of $\phi$, which, in an orthonormal frame
$\{e_{1},\dots, e_{n},{\mathbf N}\}$, is defined by: 
\begin{align*} B
&= \left(\langle\nabla_{e_{i}} e_{j}, {\mathbf N}\rangle \right)_{i,j=1,\dots,n}
= -  \left( \langle\nabla_{e_{i}} {\mathbf N}, e_{j}\rangle 
\right)_{i,j=1,\dots,n} 
\end{align*} 
we have
$$
|\nabla^{\phi}_{e_{i}}{\mathbf N}|^2 = \left<
\nabla^{\phi}_{e_{i}}{\mathbf N}, \nabla^{\phi}_{e_{i}}{\mathbf
N}\right\rangle = \sum_{j=1}^{n}\left< \nabla^{\phi}_{e_{i}}{\mathbf N},
e_{j}\right>^{2} \quad (\forall i=1,\dots,n),
$$ 
which implies that 
$$
\sum_{i=1}^{n} \langle \nabla^{\phi}_{e_{i}}{\mathbf N},
\nabla^{\phi}_{e_{i}}{\mathbf N} \rangle =  |B|^{2}.
$$

\noindent In
conclusion: 
$$ 
-\left< \tau_{2,\lambda}(\phi) , {\mathbf N}\right>  = n \left(
-\Delta H + H |B|^{2} - H \ricci({\mathbf N})+\lambda H\right) .
$$ 
\end{proof}

\begin{corollary}
An isometric immersion $\phi : M^n \to N^{n+1}(c)$ into a space form of 
constant curvature $c$ is biminimal if and only if there exists a real number $\lambda$
such that:
\begin{equation*}
\Delta H -H (n^{2}H^2 - s + n(n-2)c+\lambda) =0,
\end{equation*}
where $H$ is the mean curvature and $s$ the scalar curvature of $M^n$.
Moreover, an isometric immersion $\phi : M^2 \to N^{3}(c)$ from
a surface to a three-dimensional space form  is biminimal if
and only if: \begin{equation}\label{eq6} 
\Delta H -2H (2H^2 - G+\frac{\lambda}{2})=0,
\end{equation}
for some  $\lambda$ in $\r$.
\end{corollary}

\begin{proof}
Let $\{e_{1},\dots,e_{n}\}$ be an orthonormal frame of $M^n$ 
corresponding to 
the principal curvatures $\{k_{1},\dots,k_{n}\}$ and $B$ its second
fundamental form, then: \begin{align*}
|B|^2 &= k_{1}^2 + \cdots + k_{n}^2 = n^2 H^2 - 2 \sum_{\substack{i,j=1\\ i<j}}^{n} k_{i}k_{j} \\
&= n^2 H^2 - 2 \sum_{\substack{i,j=1\\ i<j}}^{n} (K(e_{i},e_{j}) -c)
 = n^2 H^2 - \sum_{i,j=1}^{n} K(e_{i},e_{j}) + n(n-1)c \\
&= n^2 H^2 - s + n(n-1)c ,
\end{align*}
where $K(e_{i},e_{j})$ is the sectional curvature on $M^n$ of the plane 
spanned by $e_{i}$ 
and $e_{j}$, and $s=\sum_{i,j=1}^{n} K(e_{i},e_{j})$ is the scalar
curvature of $M^n$. Since $\ricci({\mathbf N}) = nc$, the
map $\phi$ is biminimal if and only if: 
$$
\Delta H = (n^2 H^2 -s +
n(n-2)c+\lambda) H ,
$$ 
for some  $\lambda$ in $\r$.
\end{proof}

\begin{remark}
Condition~\eqref{eq6}, for free biminimal immersions, is very similar to the equation of the Willmore 
problem ($\Delta H + 2 H (H^2 - K)=0$) but the minus sign in
\eqref{eq6} rules out the existence of compact solutions when $c\leq
0$. 
\end{remark}

We shall now describe some constructions to produce examples of
biminimal immersions. Recall that a submersion
$\phi:(M,g)\to(N,h)$ between two Riemannian manifolds is
{\it horizontally homothetic} if there exists a
function $\Lambda:M\to\r$, the {\it dilation}, such
that: 
\begin{itemize} 
\item at each point $p\in M$ the differential
$d\phi_p:H_{p}\to T_{\phi(p)}N$ is a conformal map with factor
$\Lambda(p)$, i.e.
$\Lambda^2(p)g(X,Y)(p)=h(d\phi_{p}(X),d\phi_{p}(Y))(\phi(p))$ for all
$X,Y\in H_{p}=Ker_{p}(d\phi)^{\perp}$;
\item $X(\Lambda^2)=0$ for all horizontal vector fields. 
\end{itemize}   

\begin{lemma}\label{lemma1} Let $\phi : (M^n,g)\to
(N^2,h)$ be a horizontally homothetic submersion with  
$\Lambda$ and minimal fibres and  let $\gamma:
I\subset {\r} \to N^2$ be a curve parametrised by arc-length, of
signed curvature $k_{\gamma}$. Then the codimension-one submanifold
$S=\phi^{-1}(\gamma(I))\subset M$ has mean curvature
$H_{S}=\Lambda k_{\gamma}/(n-1)$.  
\end{lemma}
\begin{proof}
Let $\{ T,N\}$ be the Frenet frame of $\gamma$, 
i.e. $\nabla_{\frac{\partial}{\partial t}}^{\gamma} T =  k_{\gamma}
N$.
Choose a local orthogonal frame $\{e_{1},e_{2},e_{3},\ldots,e_{n}\}$
on $M^{n}$ such that $d\phi(e_{1}) = T\circ \phi$, $d\phi(e_{2}) =
N\circ \phi$ and $d\phi(e_{i}) = 0$, for $i=3,\dots,n$. Since $\phi$
is a horizontally homothetic submersion, we have that 
$|e_{1}|^2=|e_{2}|^2=1/\Lambda^2$ and  can
choose $\{e_{3},\ldots,e_{n}\}$ of unit length. 
The restriction to $S=\phi^{-1}(\gamma(I))\subset M$ of the
vector fields $e_{1}$ and $\{e_{3},\ldots,e_{n}\}$ give a local frame
of vector fields tangent to the submanifold  $S =
\phi^{-1}(\gamma(I))$, while the restriction of $\Lambda e_{2}$
gives a unit  vector field normal to $S$. 
Therefore the mean curvature
of $S$ is:  
$$ 
H_{S} = \frac{1}{n-1} \Lambda^3 \left< \nabla_{e_{1}} e_{1},e_{2}\right>  + 
\frac{1}{n-1} \Lambda \sum_{i=3}^{n}\left< \nabla_{e_{i}} e_{i}, e_{2}
\right> =
\frac{1}{n-1} \Lambda^3 \left< \nabla_{e_{1}} e_{1},e_{2}\right>  + 
\frac{n-2}{n-1}\Lambda (H_{\rm fibre}),
$$
where $H_{\rm fibre}$ is the mean curvature of the
fibres.  The fibres being minimal ($H_{\rm fibre}=0$), we have:
\begin{equation}\label{eq:mean-lambda}
H_{S} = \frac{1}{n-1} \Lambda^3 \langle \nabla_{e_{1}}
e_{1},e_{2}\rangle. 
\end{equation}

\noindent Moreover: 
\begin{align*} 
\Lambda^2\langle\nabla_{e_{1}} e_{1}, e_{2} \rangle &= \langle d\phi(\nabla_{e_{1}}
e_{1}), d\phi(e_{2}) \rangle = \langle d\phi(\nabla_{e_{1}} e_{1}),
N\circ\phi \rangle \\ &= \langle \nabla^{\phi}_{e_{1}} d\phi(e_{1}),
N\circ\phi\rangle =  \langle \nabla^{\phi}_{e_{1}} (T\circ\phi),
N\circ\phi\rangle, 
\end{align*} 
since $(\nabla d\phi)(e_{1},e_{1}) =
0$ for a horizontally homothetic submersion (cf.~\cite{B-W}). 
\newline Finally, 
$$
\nabla^{\phi}_{e_{1}} (T\circ\phi) = (\nabla_{d\phi(e_{1})}
T)\circ\phi =  (\nabla_{T}T)\circ\phi =k_{\gamma} N \circ\phi,
$$ 
and, taking into account \eqref{eq:mean-lambda}, $(n-1) H_{S} =
\Lambda k_{\gamma}$. 
\end{proof}

\begin{theorem}\label{teo-main1}
Let $\phi : M^3(c)\to
(N^2,h)$ be a horizontally homothetic submersion with dilation
$\Lambda$, minimal fibres and integrable horizontal
distribution, from a space form of constant sectional curvature $c$ to
a surface. Let  $\gamma: I\subset {\r} \to N^2$ be a curve parametrised
by arc-length such that the surface
$S=\phi^{-1}(\gamma(I))\subset M^3$ has constant Gaussian curvature $c$.
Then  $S=\phi^{-1}(\gamma(I))\subset M^3$  is a biminimal surface (w.r.t. $2c$)
if and only if $\gamma$ is a free biminimal curve. 
\end{theorem}

\begin{proof}
Let $\{ T,N\}$ be the Frenet frame of $\gamma$, 
i.e. $\nabla_{T} T = 
k_{\gamma}N$. Let $\tilde{\gamma}:I\to M^3$ be a horizontal lift 
of $\gamma$, so that
$\tilde{\gamma}'$ is horizontal and
$d\phi(\tilde{\gamma}')=T\circ\phi$. 
Let $\psi(t,s)=\eta_{s}(\tilde{\gamma}(t))$ be a local parametrisation
of the surface $S=\phi^{-1}(\gamma(I))\subset M^3$,
where, for a fixed $t_0\in I$, $\eta_{s}(\tilde{\gamma}(t_0))$ is a parametrisation
by arc-length of the fibre of $\phi$
through $\tilde{\gamma}(t_0)$. Then $\psi$ induces on the surface
$S$ the metric 
$$
g_{S}= \frac{1}{\Lambda^2} dt^2 + ds^2,
$$
where $\Lambda$ is the dilation of $\phi$ which, when restricted to the
surface $S$, depends only on $s$. The Laplacian on $S$
is then given by:
\begin{equation}\label{eq-laplacian-S}
\Delta=\Lambda^2 \frac{\partial^2}{\partial t^2}+
\frac{\partial^2}{\partial s^2} - \grad(\log\Lambda)
\frac{\partial}{\partial s},
\end{equation}
whilst the Gaussian curvature of $S$ reduces to:
\begin{equation}\label{eq-gauss-S}
G_{S}=\frac{\Delta\Lambda}{\Lambda}-(\grad(\log\Lambda))^2 .
\end{equation}
Now, assuming that $S$ has constant Gaussian curvature $G_{S}=c$,
from \eqref{eq6}, $S$ is biminimal (w.r.t. $2c$) in $M^3(c)$ if
and only if: 
$$ 
\Delta {H} - 2{H}(2{H}^2-c+c)= \Delta {H} -
4{H}^3=0.
$$  
By Lemma~\ref{lemma1}, $2{H}=\Lambda k_{\gamma}$, thus 
\begin{align}
2(\Delta {H} - 4{H}^3)= & \Delta (\Lambda k_{\gamma}) -
(\Lambda k_{\gamma})^3\nonumber \\ 
= &\Lambda^3 \left[k_{\gamma}''-k_{\gamma}^3 +
\frac{k_{\gamma}}{\Lambda^2}\left(G_{S}+(\grad(\log\Lambda))^2\right)
\right]\nonumber\\ 
= & \Lambda^3 \left[k_{\gamma}''-k_{\gamma}^3 +
\frac{k_{\gamma}}{\Lambda^2}\left(c+(\grad(\log\Lambda))^2\right)
\right]\label{eq-teo-hh}.
 \end{align} 
Finally, from the generalised O'Neil formula  relating the sectional
curvatures of the domain and target manifolds for a given
horizontally homothetic submersion with integrable horizontal
distribution (see, for example \cite[Corollary~11.2.3]{B-W}) we get: 
$$
\frac{1}{\Lambda^2}\left(c+(\grad(\log\Lambda))^2\right)=G_{N},
$$
which, together with \eqref{eq-teo-hh}, gives: 
$$
2(\Delta {H} - 4{H}^3)= \Lambda^3(k_{\gamma}''-k_{\gamma}^3
+k_{\gamma} G_{N}).
$$
Then the theorem follows from Corollary~\ref{cor2.1}.
\end{proof} 

When the horizontal space is not integrable, we can reformulate 
Theorem~\ref{teo-main1}, for Riemannian submersions.

\begin{theorem}\label{teo-main2}
Let $\phi : M^3(c)\to
N^2(\bar{c})$ be a Riemannian submersion
 with minimal fibres from a
space of constant sectional curvature $c$ to a surface of
constant Gaussian curvature $\bar{c}$. Let  $\gamma: I\subset {\r} \to
N^2$ be a curve parametrised by arc-length.
Then  $S=\phi^{-1}(\gamma(I))\subset M^3$  is a biminimal surface (w.r.t. $\lambda$)
if and only if $\gamma$ is a biminimal curve (w.r.t. $\lambda + \bar{c}$). 
\end{theorem}
\begin{proof}
First, from \eqref{eq-laplacian-S} and \eqref{eq-gauss-S}, since
$\Lambda=1$, we have: 
\begin{align*}
& G_{S}=0 \, , \quad \Delta=\frac{\partial^2}{\partial t^2}+ \frac{\partial^2}{\partial
s^2}. \end{align*}
Thus, taking into account Lemma~\ref{lemma1}, $S$ is biminimal if and
only if 
$$\Delta (2{H}) - (2{H})^3 - 2H\lambda =k_{\gamma}''-k_{\gamma}^3 -k\lambda=0.
$$ 

>From Corollary~\ref{cor2.1}, the latter equation is clearly biminimality (w.r.t. $\lambda + \bar{c}$) for a curve $\gamma:I\to N^2(\bar{c})$.

\end{proof}

\section{Examples of biminimal surfaces in three-dimensional space forms}

\subsection{Examples of biminimal surfaces in $\r^3$}

We apply Theorem~\ref{teo-main1} to construct examples of
biminimal surfaces in $\r^3$ with the flat metric.
\begin{enumerate}
\item
First we consider the orthogonal projection $\pi:\r^3\to\r^2$, given
by $\pi(x,y,z)=(x,y)$. The projection $\pi$ is clearly a Riemannian
submersion with minimal fibres (vertical straight lines in $\r^3$) and
integrable horizontal distribution. Thus, from Theorem~\ref{teo-main1}
a vertical cylinder with generatrix a free biminimal curve of $\r^2$ 
is a free biminimal surface. For example one can consider the cylinder on
the logarithmic spiral.  
\item
The space $\r^3\setminus\{0\}$ can be described as the warped product
$\r^3\setminus\{0\}=\r^{+}\times_{t^2}\s^2$ with the warped metric
$g=dt^2+t^2d\theta^2$, $d\theta^2$ being the canonical metric on
$\s^2$.
Then projection onto the second factor
$\pi_2:\r^{+}\times_{t^2}\s^2\to \s^2$ is a horizontally homothetic
submersion with dilation $1/t$, integrable horizontal distribution
and minimal fibres. Geometrically $\pi_2$ is the radial projection
$p\mapsto p/|p|,\;p\in\r^3\setminus\{0\}$. Applying again
Theorem~\ref{teo-main1}, we see that the cone on a
free biminimal curve on $\s^2$ is a free biminimal surface of
$\r^3$. For example if we take the parallel on $\s^2$ of latitude
$\pi/4$, which is a biharmonic curve, and thus free biminimal, we get the
standard cone of revolution in $\r^3$.
\item The following example does not seem to enter the picture of Theorem~\ref{teo-main1}. Let $\alpha: I\subset {\r} \to {\r}^{3}$
 be a space curve with 
curvature $k$ equal to its torsion $\tau$ and 
$\{T,N,B\}$ its Frenet frame. It is easy to see that the envelope $S$
of $\gamma$, parametrised by: $X(u,s) = \alpha(s) + u(B+T)$,
has mean curvature $H=k$. Thus $S$ is free biminimal if and only if: 
\begin{equation}\label{eq-envelope}
\Delta H - 4 H^3 = k'' - 4 k^3 = 0.
\end{equation}
Geometrically the curve $\gamma$ is a curve with constant slope,
i.e. there exists a vector $u\in\r^3$ such that $\langle T,u\rangle$ is constant. Then $\gamma$ can be described as an helix of the
cylinder on a plane curve $\beta$ (the orthogonal projection of
$\gamma$ onto a plane orthogonal to $u$) whose geodesic curvature is a
solution of \eqref{eq-envelope}. For example we can take $\beta$
to be the logarithmic spiral of natural equation 
$k_{\beta}=1/(\sqrt{2} s)$.
\end{enumerate}

\subsection{Examples of biminimal surfaces in $\h^3$}
\begin{enumerate}
\item Let $\h^3=\{(x,y,z)\in\r^3\,:\, z>0\}$ be the half-space model
for the hyperbolic space endowed with the metric of constant sectional
curvature $-1$ given by $g=(dx^2+dy^2+dz^2)/z^2$. Then the projection
onto
the plane at infinity defines a horizontally homothetic submersion
$\pi:\h^3\to\r^2$ with dilation $\Lambda=z$, integrable horizontal
distribution and minimal fibres (vertical lines in $\h^3$).
 Then, from Theorem~\ref{teo-main1},
a vertical cylinder with generatrix a free biminimal curve of $\r^2$  is
a biminimal surface (w.r.t. $-2$) in the hyperbolic space. For example 
the cylinder on the logarithmic spiral is free biminimal in $\r^3$ whilst
it is biminimal (w.r.t. $-2$) in $\h^3$.  
\item Let 
$\pi:\h^3\to\h^2$ defined by  $\pi(x,y,z)=(x,0,\sqrt{y^2+z^2})$. The
fibre of $\pi$ over $(x,0,r)$ is the semicircle with centre
$(x,0,r)$, radius $r$, and parallel to the coordinate
$yz$-plane. Thus the map $\pi$ has minimal fibres. 
Geometrically, this map is a projection along the geodesics of
$\h^3$ which are orthogonal to $\h^2$. This is again a  horizontally
homothetic submersion with dilation, along the fibres,
$\Lambda(s)=1/\cosh(s)$, $s$ being the arc-length parameter of the
fibre. An easy computation shows that for any curve $\gamma$
parametrised by arc-length in $\h^2$ the surface $S=\pi^{-1}(\gamma(I))$
is of constant Gaussian curvature $-1$. Then, applying
Theorem~\ref{teo-main1}, for any free biminimal curve of $\h^2$ 
$S=\pi^{-1}(\gamma(I))$ is a biminimal surface (w.r.t. $-2$) in the hyperbolic
space. 
\end{enumerate}

\subsection{Examples of biminimal surfaces in $\s^3$}
\begin{enumerate}
\item Let $p,q\in\s^3$ be two antipodal points. Then the space
$\s^3\setminus\{p,q\}$ can be described as the warped product
$\s^3\setminus\{p,q\}=(0,\pi)\times_{\sin^2(t)}\s^2$ with the warped
metric $g=dt^2+\sin^2(t)d\theta^2$, $d\theta^2$ being the canonical
metric on $\s^2$. Then the projection to the second factor
$\pi_2:\r^{+}\times_{t^2}\s^2\to \s^2$ is a horizontally homothetic
submersion with dilation $1/\sin(t)$, integrable horizontal
distribution and minimal fibres. Geometrically, $\pi_2$ is the 
projection along the longitudes onto the equatorial sphere.
Theorem~\ref{teo-main1} gives a correspondence
between free biminimal curves on $\s^2$ and biminimal surfaces (w.r.t. $2$) of
$\s^3$ given by $S=\pi_2^{-1}(\gamma(I))$.
\item This is the only example for which we use
Theorem~\ref{teo-main2}. Let $H:\s^3\to\s^2(\tfrac{1}{2})$ be the Hopf map
defined by $H(z,w)=(2z \bar{w},|z|^2-|w|^2)$, where we have
identified  $\s^3=\{(z,w)\in\c^{2}\;:\;|z|^2+|w|^2=1\}$ and 
$\s^2(\tfrac{1}{2})=\{(z,t)\in\c\times\r\;:\;|z|^2+t^2=\tfrac{1}{4}\}$. The Hopf map
is a Riemannian submersion with minimal fibres (great circles in
$\s^3$).  Thus, from Theorem~\ref{teo-main2}, we see that 
a Hopf cylinder $H^{-1}(\gamma(I))$ is a biminimal surface (w.r.t. $\lambda$) of $\s^3$ if and only if the curve $\gamma$ is a biminimal curve (w.r.t. $\lambda +4$) of $\s^2(\tfrac{1}{2})$.
\end{enumerate}

\section{Examples of biminimal surfaces in  
Thurston's three-dimensional geometries}

Of Thurston's eight geometries (cf.~\cite[Section~10.2]{B-W}), three have 
constant sectional curvature, $\r^3$, $\s^3$ and $\h^3$, and
contain biminimal surfaces as described in the previous section, 
two are  Riemannian products, $\s^2 \times \r$ and $\h^2 \times \r$,
and will be our first class of examples, two are line bundles, 
over $\r^2$ for ${\mathcal H}_3$ and over $\r^2_{+}$ for $\sl$, and one,
$\sol$, does not allow Riemannian submersion or horizontally 
homothetic maps with minimal fibres to a surface, even locally,
and therefore does not fit our framework.

\subsection{Biminimal surfaces of $\s^2 \times \r$ and $\h^2 \times \r$}
In both cases, consider the Riemannian submersion with totally geodesic fibres, given 
by the projection onto the first factor, $\pi : N^2 \times \r \to N^2$. Given a curve 
$\gamma : I \subset {\r} \to N^2$ parametrised by 
arc length, take its Frenet frame $\{ T,N\}$, and 
consider $\{
e_{1},e_{2}\} \in T(N^2\times\r)$ its horizontal lift. 
The unit vertical vector $e_{3}$ completes $\{ e_{1},e_{2}\}$
 into an orthonormal
frame of $T(N^2 \times \r)$, such that $\{ e_{1},e_{3}\}$ 
is a basis of $TS$, for $S=\pi^{-1}(\gamma(I))$, with $e_{2}$ the
normal to the surface. Then, from Lemma~\ref{lemma1} the 
mean curvature of $S$ is $H= \tfrac{k}{2}$, where $k$ is the 
signed curvature of
$\gamma$, and, from Proposition~\ref{pro-B-ricci}, $S$ is 
biminimal (w.r.t. $\lambda$) in $N^2 \times \r$ if:
$$ \Delta H = (|B|^2 - \ricci(e_{2}) +\lambda) H.$$
With respect to the frame $\{ e_{1},e_{3}\}$  the 
matrix associated to the second fundamental form of $S$ is:
$$
B=\begin{pmatrix}
k & 0 \\
0 & 0
\end{pmatrix}.
$$
Besides, 
$$\ricci^{N^2 \times {\tiny \r}}(e_{2}) = \ricci^{N^2}(e_{2}) =
\left\lbrace
\begin{array}{l}
+ 1 \quad \mbox{if } \, N^2 = \s^2 \\
- 1 \quad \mbox{if } \, N^2 = \h^2.
\end{array}
\right.
$$
In both cases, using Equation~\eqref{eq-laplacian-S}, 
$\Delta H = \Delta (\tfrac{k}{2}) = \frac{1}{2} k''$, so $S$ is
 biminimal in $N^2 \times \r$ if:
$$ k'' = k^3 - k + \lambda k\quad \mbox{if } \, N^2 = \s^2 $$
and 
$$ k'' = k^3 + k + \lambda k\quad \mbox{if } \, N^2 = \h^2 .$$
Now comparing with \eqref{eq1}, we have the following 
\begin{proposition} The cylinder $S=\pi^{-1}(\gamma(I))$ 
is a biminimal surface (w.r.t. $\lambda$) in $N^2 \times \r$  if and only 
if $\gamma$ is a biminimal curve (w.r.t. $\lambda$) on $N^2$ ($\s^2$ or $\h^2$) .
\end{proposition}

\subsection{Biminimal surfaces of the Heisenberg space}
The three-dimensional Heisenberg space ${\mathcal H}_3$ is the
two-step nilpotent Lie group standardly represented in $GL_3(\r)$ by 
$$ 
\left[
\begin{array}{ccc} 
1&x&z\\ 0&1&y\\ 0&0&1 
\end{array}
\right] 
$$ 
with $x,y,z\in\r$. 
Endowed with the left-invariant metric
\begin{equation}\label{metric} 
g = dx^{2} + dy^{2} + \left(dz - xdy\right)^{2}, 
\end{equation} $({\mathcal H}_3,g)$ has a rich geometric
structure, reflected by the fact that its group of isometries 
is of dimension $4$, the maximal possible dimension 
for a metric of non-constant curvature on a three-manifold. 
Also, from the algebraic point of view, this is a two-step 
nilpotent Lie group, i.e. ``almost Abelian''. 
An orthonormal basis of left-invariant vector fields is given, 
with respect to the coordinates vector fields, by

\begin{equation}\label{inv-fields} 
E_{1}= \frac{\partial}{\partial x};
\quad 
E_{2}=\frac{\partial}{\partial y}+ x \frac{\partial}{\partial z};
\quad 
E_{3}=\frac{\partial}{\partial z}.
\end{equation} 

Let now $\pi:{\mathcal H}_3\to\r^2$ be the projection 
$(x,y,z)\mapsto (x,y)$. At a point
$p=(x,y,z)\in {\mathcal H}_3$ the vertical space of the 
submersion $\pi$ is
$V_p=\rm{Ker}(d\pi_p)=\rm{span}(E_3)$ and the 
horizontal space is 
$H_p=\rm{span}(E_1,E_2)$. An easy computation shows 
that $\pi$ is a Riemannian submersion with minimal 
fibres. 

\noindent Take a curve $\gamma(t)=(x(t),y(t))$ in $\r^2$, 
parametrized by arc length, with signed curvature $k$, and consider the flat cylinder
 $S=\pi^{-1}(\gamma(I))$ in ${\mathcal H}_3$.
Since the left invariant vector fields are orthonormal, 
the vector fields
$$
e_1=x'E_1+y'E_2;\quad e_2=E_3
$$  
give an orthonormal frame tangent to $S$ and 
$$
N=-y'E_1+x'E_2
$$
is a unit normal vector field of $S$ in ${\mathcal H}_3$. The second fundamental form of $S$ is
$$
B=\begin{pmatrix}
k & -\frac{1}{2} \\
-\frac{1}{2} & 0
\end{pmatrix}.
$$
Clearly $H=\trace(B)/2=k/2$, $|B|^2=k^2+1/2$ and
a direct computation shows that 
$\ricci(N)=-\tfrac{1}{2}$. Thus, from \eqref{eq5}, $S$ is biminimal w.r.t. $\lambda$ if and only if
$$
\Delta H =(|B|^2-\ricci(N)+\lambda)H
$$
or equivalently
$$
k''=(k^2+1/2+1/2+\lambda)k=k^3+k(1+ \lambda) .
$$
Finally, taking \eqref{eq1} into account, we have:

\begin{proposition}
The flat cylinder 
$S=\pi^{-1}(\gamma(I))\subset {\mathcal H}_3$ is a biminimal surface (w.r.t. $\lambda$) of 
${\mathcal H}_3$ if and only if $\gamma$ is a biminimal curve (w.r.t. $\lambda + 1$) of $\r^2$. 
\end{proposition}

\subsection{Biminimal surfaces of $\sl$}

Following \cite[page 301]{B-W} we identify $\sl$ with: 
$$
\r^3_{+}=\{(x,y,z)\in\r^3\,:\,z>0\}
$$
endowed with the metric:

\begin{equation}\label{metric-sl}
ds^2=\left(dx+\frac{dy}{z}\right)^2+\frac{dy^2+dz^2}{z^2}.
\end{equation}

Then the projection $\pi:\sl\to \r^2_{+}$ defined by 
$(x,y,z)\mapsto(y,z)$ is a submersion and if we denote, as usual, 
by $\h^2$
the space
$\r^2_{+}$ with the hyperbolic metric $\frac{dy^2+dz^2}{z^2}$,
the submersion $\pi:\sl\to \h^2$ becomes a Riemannian submersion 
with minimal fibres.
The vertical space at a point $p=(x,y,z)\in\sl$ is 
$V_p=\rm{Ker}(d\pi_p)=\rm{span}(E_1)$ and the horizontal space at
$p$ is $H_p=\rm{span}(E_2,E_3)$, where 

\begin{equation}\label{ort-fields} 
E_{1}= \frac{\partial}{\partial x};
\quad 
E_{2}=z\frac{\partial}{\partial y}- \frac{\partial}{\partial x};
\quad 
E_{3}=z\frac{\partial}{\partial z}
\end{equation} 
give an orthonormal frame on $\sl$ with respect to the metric 
\eqref{metric-sl}. 

\noindent Now, given a curve $\gamma(t)=(y(t),z(t))$ on $\h^2$, 
parametrized by arc length, and the flat cylinder $S=\pi^{-1}(\gamma(I))$ in $\sl$,
as $E_{1},E_{2}$ and $E_{3}$ are orthonormal, the vector fields
\begin{equation}\label{orto-frame}
e_1=\frac{y'}{z}E_2+\frac{z'}{z}E_3;\quad e_2=E_1
\end{equation} 
give an orthonormal frame tangent to $S$ and 
$$
N=-\frac{z'}{z}E_2+\frac{y'}{z}E_3
$$
is a unit normal vector field of $S$ in $\sl$.

\noindent With  calculations similar to those of the previous 
example, we find
that, with respect to the orthonormal frame~\eqref{orto-frame}:
$$
B=\begin{pmatrix}
k & \frac{1}{2} \\
\frac{1}{2} & 0
\end{pmatrix},\quad \ricci(N)=-\frac{3}{2}.
$$

\noindent Thus:

\begin{proposition}
The flat cylinder $S=\pi^{-1}(\gamma(I))\subset \sl$ is a biminimal surface (w.r.t. $\lambda$)
of $\sl$ if and only if $\gamma$ is a biminimal curve (w.r.t. $\lambda + 1$) of $\h^2$. 
\end{proposition}

\begin{remark}
These links between biminimal cylinders and biminimal curves are very similar to 
those described by U.~Pinkall~\cite{Pink85} between Willmore Hopf tori of $\s^3$ and
elastic curves on $\s^2$.
\end{remark}

\end{document}